\documentclass{amsart}
\usepackage[all]{xy}
\usepackage{graphicx}
\usepackage{psfrag}
\usepackage{amsthm}
\usepackage{amscd}
\usepackage{amsfonts}
\usepackage{amstext}
\usepackage{amssymb}
\usepackage{amsmath}
\usepackage{amsxtra}
\usepackage{plain}
\usepackage[all]{xy}

\newtheorem{prop}{\bf Proposition}[section]

\newtheorem{lem}[prop]{{\bf Lemma}}
\newtheorem{thm}[prop]{{\bf Theorem}}

\numberwithin{equation}{section}

\newenvironment{defn}{{\bf Definition: }}{ }

\newenvironment{pf}{{\bf proof }}{\qed\endtrivlist}

\newcommand{\Z}{\mathbb{Z} }

\newcommand{\C}{\mathbb{C} }
\newcommand{\R}{\mathbb{R} }

\newcommand{\ind}{\operatorname{ind}}
\newcommand{\Id}{\operatorname{Id}}

\newcommand{\U}{\mathcal{U}}

\newcommand{\End}{\operatorname{End}}

\newcommand{\im}{\operatorname{Im}}
\newcommand{\Dom}{\operatorname{Dom}}
\newcommand{\spec}{\operatorname{spec}}

\long\def\symbolfootnote[#1]#2{\begingroup%
\def\thefootnote{\fnsymbol{footnote}}\footnote[#1]{#2}\endgroup}

\begin{document}

\title{A note on some classical results of Gromov-Lawson}
%\title{Index theory and partitioning by enlargeable hypersurfaces}
\author{Mostafa Esfahani Zadeh}
\address{Mostafa Esfahani Zadeh.\newline
 Mathematisches Institute, 
Georg-August-Universit\"{a}t,\newline
G\"{o}ttingen, Germany \and \newline 
Institute for Advanced Studies in Basic Siences(IASBS),\newline 
Zanjan-Iran}
\symbolfootnote[0]{\emph{2000 Mathematics Subject Classification}. 58J22 
(19K56 46L80 53C21 53C27).} 
\symbolfootnote[0]{\textit{Key words and phrases}. Higer index theory, 
enlargeablity, Dirac operators.}
\email{zadeh@uni-math.gwdg.de}

\begin{abstract}
In this short note we show how the higher index theory can be used to prove 
results concerning the non-existence of complete riemannian metric with uniformly 
positive scalar curvature at infinity. 
By improving some classical results due to M. Gromov and B. 
Lawson we show the efficiency 
of these methods in dealing with such non-existence theorems. 
\end{abstract}
\maketitle
\section{Introduction}
Let $(M,g)$ be an oriented complete non-compact manifold partitioned by a 
compact hypersurface $N$ into two part $M_+\text{ and } M_-$ with 
$M_+\cap M_-=\emptyset$ and $\overline M_+\cap\overline M_-=N$.
We assume that the positive unite normal to $N$ points out from $M_-$ to $M_+$.  
Let $W$ be a Clifford bundle on $M$ which is at the same time a Hilbert 
$A$-module bundle. 
We assume that this bundle is equipped with a connection which is compatible 
to the Clifford action of $TM$ and denote the corresponding $A$-linear Dirac 
type operator by $D$. As an example let $M$ be a spin manifold with spin 
bundle $S$ and $V$ be a Hilbert $A$-module bundle on $M$ with a hermitian connection. 
Then the spin Dirac operator twisted by $V$ is an example of such operator 
which acts on the smooth sections of $W=S\otimes V$.

Let $U=(D-i)(D+i)^{-1}$ be the Cayley transform of $D$ which is a $A$-linear 
bounded operator on $H=L^2(M,W)$. 

Let $\phi_+$ be a smooth function on $M$ which coincide with the characteristic 
function of $M_+$ outside a compact set and put $\phi_-:=1-\phi_+$.
It turns out that 
the operators $U_+=\phi_-+\phi_+U$ is $A$-Fredholm in the sense of Fomenko-Mischenko. 
The Fomenko-Mischenko index $\ind(U_+)$ does not depend on $\phi_+$ but on the 
cobordism class of the partitioning manifold $N$. This index is denoted by $\ind(D,N)$. 
A basic property of this index is the following c.f. \cite[theorem 2.4]{Zadeh7}:

{\it If $M$ is spin and $W=S\otimes V$, where $V$ is a flat Hilbert $A$-module bundle then 
$\ind(D,N)=0\in K_0(A)$ provided that the scalar curvature of $g$ is uniformly 
positive.}

The Clifford action of $i\vec n$ provides a $\Z_2$ grading for $W_{|N}$ and makes  
$W_{|N}$ a graded Clifford bundle on $N$. 
Let $D_N$ denote the associated Dirac type operator which acts on smooth 
sections of $W_{|N}$. It is a $A$-linear elliptic operator and has 
the Fomenko-Mischenko index $\ind D_N\in K_0(A)$. The following equality 
generalizing a result 
due to J. Roe and N. Higson \cite{Roe-partitioning, Higson-note} is 
proved in \cite{Zadeh7} 
\begin{equation}\label{asli}
\ind D_N=\ind(D,N)~.
\end{equation}
This equality has been used in \cite{Zadeh7} to prove the following 

\emph{If a complete spin manifold $(M,g)$ is partitioned by an enlargeable hypersurface 
$N$ and if there is a smooth map $\phi:M\to N$ whose restriction to $N$ is of 
non-zero degree, 
then the scalar curvature of $g$ can not be uniformly positive.}\\
In this short note we improve this result by showing that under the same conditions 
the scalar curvature of $g$ can not be uniformly positive even outside a 
compact subset of $M$. 
We will use this stronger result to improve some classical results due to M. Gromov and 
B. Lawson. The section \ref{dovm} deals with some analytical aspect of regular 
operators 
on Hilbert Modules. Here we construct the wave operator for the spin Dirac 
operators twisted by flat Hilbert module bundles and prove its unit speed property. 
In the forthcoming section we use these results to prove the vanishing 
theorem \ref{strvan} and its implication in improving some classical results 
due to M.Gromov and B.Lawson.

Here we give the definition of the enlargeability as it is introduced by 
Gromov and Lawson in \cite{GrLa3}.
\begin{defn}\label{enlar}
Let $N$ be a closed oriented manifold of dimension $n$ with  
a fixed riemannian metric $g$. The manifold $N$ is enlargeable if for each 
real number $\epsilon>0$ there is a riemannian spin cover $(\tilde N,\tilde g)$,  
with lifted metric, and a smooth map $f:\tilde N\to S^n$ such that: 
the function $f$ is constant outside a compact subset $K$ of $\tilde N$; 
the degree of $f$ is non-zero;  and
the map $f:(\tilde N,~\tilde g)\to(S^n,g_0)$ is $\epsilon$-contracting, where
$g_0$  is the standard metric on $S^n$. Being $\epsilon$-contracting means that 
$\|T_xf\|\leq\epsilon$ for each $x\in \tilde N$, where 
$T_xf\colon T_x\tilde N\to T_{f(x)}S^n$. 
The manifold $N$ is said to be area-enlargeable if the function 
$f$ is $\epsilon$-area 
contracting. This means $\|\Lambda^2T_xf\|\leq\epsilon$ for each 
$x\in \tilde N$, where 
$\Lambda^2T_xf\colon \Lambda^2T_x\tilde N\to \Lambda^2T_{f(x)}S^n$. 
\end{defn}
An enlargeable manifold do 
not admit riemannian metrics with positive scalar curvature. The relevance of 
this theorem will be clear by noticing that enlargeability depend only on the homotopy 
type of $M$ not on differential structures used in its definition while the existence 
of a metric with positive scalar curvature does depend in general on the 
underlying  differential structure, c.f. \cite{LaMi}. 

{\bf Acknowledgment:} The author would like to thanks J. Roe and T. Schick for 
helpful hints. 
\section{Functional calculus of the Dirac operator}\label{dovm}

Let $\mathbb H$ be a Hilbert $A$-module 
with $A$ a $C^*$-algebra and let $T$ be an $A$-linear map which is defined on 
a dense subspace $\Dom(T)$ of $\mathbb H$. 
The graph of $T$ is the following subset of $\mathbb H\oplus \mathbb H$
\[graph (T):=\{(u,T(u))|u\in \Dom(T)\}\]
The closure of this graph with respect to norm topology in 
$\mathbb H\oplus \mathbb H$ turns out to be 
the graph of an operator $\bar T$ which is called {the closure} of $T$. The domain 
$\Dom(\bar T)$ of this closure is a closed $A$-subspace of $\mathbb H$.
The adjoint of $T$ is the closed operator 
$T^*$ such that $\langle Tu,v\rangle=\langle u,T^*v\rangle$ for $u\in\Dom(T)$ 
and $v\in\Dom(T^*)$. 
Let $\Dom(T)=\Dom(T^*)$. Following \cite{Dan} $T$ is said to be self-adjoint 
if $T^*=T$ and 
it is said to be normal if 
\[\langle Tu,Tv\rangle=\langle T^*u,T^*v\rangle,\text{ for }u,v\in\Dom(T).\]
The operator $T$ is called regular if there is a bounded adjointable operator 
$P\in\mathcal L_A(\mathbb H\oplus \mathbb H,\mathbb H)$ with $\im P=\Dom\bar T$ 
(c.f proposition 5 of \cite{Dan}).
Regularity and self adjointness of the operator $T$ make it possible to associate to 
each continuous (not necessarily bounded) function $f$ on $\spec(T)$ a closed $A$-linear 
operator $f(T)$ on $\mathbb H$. This correspondence define a continuous 
functional calculus for $T$. 
This construction has been worked out in \cite{Baaj}. Here we follow the more 
geometric approach 
of \cite{Dan}. 
The key observation of \cite{Dan} is that the transformation 
$T\to Q(T)=T(1+T^*T)^{-1/2}$ provides a $*$-preserving bijection between 
the set of all normal regular $A$-linear operator $\mathcal R_A(\mathbb H)$ 
and the following set 
\[\mathcal V(\mathbb H)=\{Q\in \End_A(\mathbb H)|\|Q\|\leq1\text{ and }\im(1-Q^*Q)
\text{ and }\im(1-QQ^*)\text{ are dense}\}.\]
Let $\mathbb D$ denote the open unit disc in $\C$.
Each function $g\in C_0(\C)$ determines a function  $\bar g\in C_0(\mathbb D)$ 
by the following  
relation 
\[\bar g(z(1+|z|^2)^{-1/2})=g(z)(1+|g(z)|^2)^{-1/2}.\]
Clearly $\|\bar g\|\leq1$, so the bounded operator $\bar g(Q(T))$ belongs again to 
$\mathcal V_A(\mathbb H)$ and corresponds to a unique operator in 
$\mathcal R_A(\mathbb H)$ 
which is defined to be $g(T)$. 
From this construction it is clear that if $g=g'$ on a 
closed subset containing the spectrum $\spec(T)$ then $g(T)=g'(T)$.  
It is clear also that the corresponding $g\to g(T)$ provides a $C^*$-representation 
$\phi:C_0(X)\to\End_A(\mathbb H)$. 

The regularity of $T$ implies $C_0(X)(T)\mathbb H=\mathbb H$ 
(see corollary 14 and theorem 15 of \cite{Dan}), so given any $u\in \mathbb H$ there is 
$g\in C_0(X)$ and $v\in \mathbb H$ with $u=g(T)v$. If $h$ is a bounded continuous 
function on 
$X$  then $hg\in C_0(X)$ and one defines $h(T)u:=(hg)(T)v$. 
It is easy to verify that $h(T)$ is well defined and that it is 
bounded with $\|h(T)\|\leq\|h\|$, 
where $\|h\|=\sup\{h(\lambda)|\lambda\in spec(T)\}$.

For an unbounded continuous function $h$, the set 
\[\mathcal C(h):=\{g\in C_0(X)|hg\in C_0(X)\}\]
is a $*$-subalgebra of $C_0(X)$. The above argument can be used to define the unbounded 
operator $h(T)$ with domain $\phi(\mathcal C(h))(\mathbb H)$. 
It is clear that $\Dom f(T)\subset\Dom g(T)$ if $|g|\leq|f|$ on $\spec(T)$. 
As an example if $g(z)=z$ then it turns out that $g(T)=T$. The following theorem 
summarizes some  properties of this functional calculus which are relevant to our 
purposes. These properties are straightforward consequences of the above 
discussion.

\begin{thm}\label{funccal}
Let $T$ be a densly defined regular normal $A$-linear operator on Hilbert 
$A$-module $\mathbb H$ and let $X\subset \C$ be a closed subset containing 
the spectrum $spec(T)$. 
Then to each continuous function $f\in C(X)$ one can correspond the regular 
normal operator 
$f(T)$ on $\mathbb H$ satisfying $f(T)^*=\overline{f(T^*)}$ 
such that 
\begin{enumerate}
 \item for $f$, $g$ in $C(X)$ the operator $(f+g)(T)$ is the closure of 
$f(T)+g(T)$, $(fg)(T)$ is 
the closure of $f(T)g(T)$ and  $f\circ g(T)=f(g(T))$.
 \item if $T$ is bounded then $f\to f(T)$ coincides to the functional calculus 
in the $C^*$-algebra of bounded operators on $\mathbb H$.
 \item let $\{f_k\}_k$ be a sequence of continuous functions on $X$ which is 
dominated be a 
continuous function $F$, i.e. $|f_k(z)|\leq|F(z)|$ for all $z\in X$. 
If $f_k\to f$ uniformly on compact 
subsets of $X$ then $f_k(T)(u)\to f(T)(u)$ for $u\in\Dom(F)$. 
If $f_k$'s are uniformly bounded and the convergence to $f$ is uniform 
then $f_k(T)\to f(T)$ in norm topology.
\end{enumerate}
\end{thm}
As a special case of this theorem, if $T$ is self adjoint then  
$f(t,T)=e^{itT}$,  for $t\in\R$, is a one parameter family of unitary 
operators on $\mathbb H$ 
satisfying  $e^{i(t+s)T}=e^{itT}e^{isT}$. 
Moreover $Te^{itT}=e^{itT}T$ which implies that $e^{itT}$ provides a unitary bijection 
between $\Dom T$ and $\im T$. 
By the third part of the previous theorem for $u\in Dom(T)$ we have 
\[\left(\frac{d}{dt}\right)_{t=0}e^{itT}u=\lim_{t\to0}\left(\frac{e^{itx}-1}{t}(T)
\right)(u)=iT(u).\]
So, for a self adjoint regular operator $T$ the $t$-parameterized family of 
unitary operators $e^{itT}$ satisfies the wave equation and the initial condition
\begin{gather}
(\frac{d}{dt}-iT)e^{itT}(u)=0~~;\text{ for } u\in\Dom(T),\label{wave}\\
 \lim_{t\to0}e^{itT}u=u.\label{inicon}
\end{gather}
The wave equation \eqref{wave} implies  the following vanishing result for $u\in\Dom(T)$
\begin{align*}
 \frac{d}{dt}\langle e^{itT}(u),e^{itT}(u)\rangle
&=\langle iTe^{itT}(u),e^{itT}(u)\rangle+\langle e^{itT}(u),iTe^{itT}(u)\rangle\\
&=0.
\end{align*}
This conservation law and initial condition \eqref{inicon} prove the uniqueness of heat 
operator with properties \eqref{wave} and \eqref{inicon}.
If $f$ is a smooth function in the Schwartz space $\mathcal S(\R)$ then 
the Fourier transform 
$\hat f$ is in $\mathcal S(\R)$ and 
\begin{equation*}
 f(x)=\frac{1}{2\pi}\int_{-\infty}^\infty\hat f(s)e^{isx}\,ds.
\end{equation*}
By applying the above theorem we get the following formula  
\begin{equation}\label{bradar}
f(T)(u)=\frac{1}{2\pi}\int_{-\infty}^\infty\hat f(s)e^{isT}(u)\,ds~~,\
\text{ for }u\in\Dom(T).
\end{equation}
%\Red{An element $f\in\mathcal S^0(\R)$ is a tamed distribution and 
%its Fourier transform is a tamed distribution too. So the relation 
%\eqref{bradar} can be used as the 
%definition of $f(T)$. } 
The Hilbert module that we shall study in the sequel is the space of the 
$L^2$-sections of Hilbert module bundles over complete manifolds. 
Let $(M,g)$ be a complete riemannian manifold  
and let $W$ be a Clifford Hilbert $A$-modules bundle over $M$, where $A$
is a complex $C^*$-algebra. For $\sigma$ and $\eta$ two compactly 
supported  smooth sections of $W$ put 

\[\langle\sigma,\eta\rangle=\int_M\langle\sigma(x),\eta(x)\rangle\,d\mu_g(x)\in A.\] 
It is easy to show that 
$|\sigma|=\|\langle\sigma,\sigma\rangle\|^{1/2}$ is a norm on $C_c(M,W)$. 
The completion of $C_c^\infty(M,W)$ with respect to this norm is the Hilbert $A$-module 
$\mathbb H=L^2(M,W)$.
Let $D$ be an $A$-linear Dirac type operator acting on the compactly supported 
smooth sections of $W$ which form a dense sub-space of $\mathbb H$. We 
recall that $D$ is 
formally self adjoint, i.e. $\langle D\sigma,\eta\rangle=\langle\sigma,D\eta\rangle$ 
for $\sigma$ and $\eta$ as in above. Moeover $D$ is a regular operator 
(see e.g. \cite[lemma 2.1]{Zadeh7}), so one can apply the theorem \ref{funccal} 
to $D$ and define the 
bounded operator $f(D)$ on $L^2(M,W)$ for each bounded continuous function $f$ on $\R$. 
In particular we can define the wave operator $e^{itD}$. 
In the following lemma we describe a context in which the wave operator 
has finite propagation speed. 
\begin{lem}\label{fpsp}
Let $W=S\otimes V$ be the spin bundle $S$ twisted by the flat Hilbert 
$A$-module bundle $V$. The wave operator $e^{itD}$ has unite propagation speed. 
\end{lem}
\begin{pf}
To prove the assertion we give an another construction for the wave operator 
which satisfies the 
unite propagation speed. Then the uniqueness of the wave operator implies 
the desired assertion.  
In this proof we denote by $V_0$ the fiber of $V$ which is a Hilbert $A$-module.
Let $\{U_\alpha,\phi_\alpha^S\otimes\phi_\alpha^V\}_\alpha$ be a trivializing atlas 
for $M$ such that $W_{|U_\alpha}\simeq U_\alpha\times S\otimes_{\C}V_0$. 
Since $V$ is flat, we can assume that the transition functions  
$\phi_{\alpha\beta}^V:U_\alpha\cap U_\beta\to \End_A(V_0)$ are locally constant. 
Let $D'$ denote the spin Dirac operator. Then $D=D'\otimes \Id_{V_0}$  on 
smooth sections $\xi=\sum_\alpha\xi^\alpha$ of $W$ where $\xi_\alpha$ is 
supported in $U_\alpha$ and $\xi^\alpha=s(x)\otimes v$ for a fixed $v$ in $V_0$. 
From the unite propagation formula for $D'$, for $|t|$ sufficiently small the 
section $\xi_t^\alpha(x):=e^{itD'}s(x)\otimes v$ is supported in $U_\alpha$ too. 
Moreover it is the unique 
solution of the following wave equation with the given initial condition $\xi$
\[(\frac{d}{dt}-iD'\otimes\Id_{V_0})\xi_t^\alpha(x)=0.\]
Since $\phi_{\alpha\beta}^V$ is constant, the transition of solution $\xi_t^\alpha$ to 
another chart is the solution of the wave equation in that chart with transitioned 
initial condition. Therefore these local solutions of local wave equations actually 
paste togeather to define a global solution $\xi_t$ of the wave equation for $D$ 
for sufficiently small values of $t$.
We can use $\xi_t$ as he initial condition and repeat the above procedure 
to define the solution 
of the wave equation beyond $t$. This way we get a solution which is defined for $t\in\R$. 
We define $e^{itD}\xi_0$ to be $\xi_t$.  From this construction it 
is clear that the wave operator $e^{itD}$ has unite propagation speed. 
\end{pf}
For next uses we rewrite the relation \eqref{bradar} with the 
Dirac operator $D$
\begin{equation}\label{bradar1}
f(D)(u)=\frac{1}{2\pi}\int_{-\infty}^\infty\hat f(s)e^{isD}(u)\,ds~~,\
\text{ for }u\in\Dom(D).
\end{equation}

\section{Vanishing theorem and its implications}
As mentioned in above the higher index $\ind(D,N)$ vanishes if the 
scalar curvature of the underlying riemannian metric $g$ is uniformly positive. 
In fact this 
index vanishes even if the scalar curvature is uniformly positive outside a 
compact subset of $M$. 
More precisely we prove the following theorem 
\begin{thm}\label{strvan}
With above notation if the scalar curvature of $g$ is uniformly positive at 
infinity and if $W=S\otimes V$, where $V$ is a flat Hilbert $A$-module bundle, 
then $\ind(D,N)=0\in K_0(A)$.  
\end{thm}
\begin{pf}
Let $U_0$ and $U_1$ be disjoint open subsets of $M$ such that the closure of $U_0$ 
is compact and $M=\bar U_0\cup U_1$.  
Moreover we assume that the scalar curvature $\kappa$ of $g$ is uniformly 
positive in $U_1(2r)$, e.g $\kappa>4\kappa_0$. 
Here $U_1(2r)$ consists of all point suited within distance $2r\geq0$ from $U_1$ 
and $r$ is a sufficiently big number that will be determind in below.
By multiplying the metric $g$ with a sufficiently small positive number we 
can and will assume that the constant $\kappa_0$ is arbitrary big.  
Let $\phi_0$, $\phi_1$ and $\phi_r$ be respectively the characteristic functions of 
$\bar U_0$, $U_1$ and $U_1(r)$ in $M$. 
If $\phi$ is a function on $M$ which is locally constant outside a compact subset 
then $[(D+i)^{-1},\phi]$ is compact, c.f. \cite[lemma 2.2]{Zadeh7}. Therefore
\begin{align*}
U_+&=\Id-2i\phi_+\sum_{i,j=0}^1\phi_i(D+i)^{-1}\phi_j\\
   &\sim \Id-2i\phi_+\phi_1(D+i)^{-1}\phi_r
\end{align*}
The function $(x+1)^{-1}$ can be uniformly approximated by compactly supported smooth  
functions, and hence by smooth functions with compactly supported Fourier 
transform. Let $h$ be such a function whose Fourier transform $\hat h$ is 
supported in $[-r,r]$ (the $r$ at the beginning of the proof is determined here). 
If $h$ is sufficiently close to $(x+i)^{-1}$ in sup-norm, 
then $\Id-2i\phi_+\phi_1h(D)\phi_r$ being close to 
$\Id-2i\phi_+\phi_1(D+i)^{-1}\phi_r$ in operator norm (c.f. theorem \ref{funccal}), 
is an $A$-Fredholm operator with the same index. 
Therefore we need to prove the vanishing of the index of the following operator 
\begin{equation}\label{jagozari}
\Id-2i\phi_+\phi_1h(D)\phi_r~.
\end{equation}
Let $\sigma$ be a smooth section of $W=S\otimes V$ supported in $U_1(r)$. 
Since $V$ is flat the following generalized Lichnerowicz formula holds 
with respect to the $A$-valued $L^2$-inner product, c.f \cite[page 199]{J.Ros1}
\[D^2=\nabla^*\nabla+\frac{\kappa}{4}\] 
which implies  
\begin{align*}
\langle D\sigma,D\sigma\rangle&=\langle D^2\sigma,\sigma\rangle\\
&=\langle\nabla\sigma,\nabla\sigma\rangle+\langle\frac{\kappa}{4}\sigma,\sigma\rangle\\ 
&\geq \kappa_0\|\sigma\|^2
\end{align*}
In the last inequality we have used the fact that for $a$ and $b$ in a 
$C^*$-algebra, $a+b\geq b$ and $\|a+b\|\geq\|b\|$ provided that $a$ and $b$ are positive 
and self adjoint. Consider the restriction of the Dirac operator $D$ to the 
Hilbert $A$-module $\mathcal H:=L^2(U_1(2r),W)$ and denote it by 
$D_{2r}$. This is an unbounded 
operator acting on smooth sections compactly supported in $U_1(2r)$. This operator is 
symmetric and satisfies the above positivity condition. In fact it has a 
self-adjoint regular extension to $\mathcal H$ as we are going to show. 
We recall here the notation of the proof of the lemma \ref{fpsp}. We assume the 
trivializing charts $(\U_\alpha,\phi_\alpha^S)$ and $(U_\alpha,\phi_\alpha^V)$ 
for vector bundle $S$ and $V$ over $M$ such that the transition 
functions $\phi_{\alpha\beta}=\phi_\alpha^V\circ(\phi_\beta^V)^{-1}$ from 
$U_\alpha\cap U_\beta$ into $\End_A(V_0)$ are constant. 
Since the twisting bundle $V$ is flat, the Hilbert $A$-module 
$\mathcal H$ is generated by elements 
$s\otimes v$ where $s$ is a smooth section 
of the spin bundle $S\to U_1(2r)$ supported in one of $U_\alpha$'s and $v$ is a 
constant element of $V_0$ (i.e. a flat section of $V_{|U_\alpha}$). 
On these sections the operator $D_{2r}$ 
takes the form $D'\otimes \Id$, where $D'$ denotes the spin Dirac operator 
acting on the smooth sections of $S$ which are compactly supported in $U_1(2r)$. 
With this domain, $D'$ is a symmetric operator on $L^2(U_1(2r),S)$ satisfying the 
following positivity relation in $\R$
\[\langle D's,D's\rangle\geq\kappa_0\|s\|^2~.\] 
The Friedrichs' extension theorem provides a self adjoint extension 
$\bar D'$ of $D'$ to $L^2(U_1(2r),S)$ satisfying still the above positivity condition. 
We recall two fact from the construction of the Friedrichs' extension. 
If $s$ is compactly supported in 
$U_\alpha$ then $\bar D's$ is compactly supported in $U_\alpha$ too. 
Moreover if $s$ belongs to $\Dom(\bar D')$  then $\phi\,s\in\Dom(\bar D')$ 
for a smooth compactly supported function $\phi$.  
Now define the operator $\bar D_{2r}$ as follows. Its domain consists of  
sum of sections $s\otimes v$ where $s$ is supported in $U_\alpha$ and belongs 
to $\Dom(\bar D')$ and $v$ is an element of $V$. On theses sections we define   
$\bar D_{2r}(s\otimes v)=\bar D's\otimes v$. 
This is a self adjoint operator on $\mathcal H$ satisfying the following relation 
\[\langle \bar D_{2r}\sigma,\bar D_{2r}\sigma\rangle\geq\kappa_0\|\sigma\|^2~;
~\text{ for }\sigma\in\Dom(\bar D_{2r})\subset\mathcal H\]
Since $\bar D'$ is self adjoint $\im(\bar D'+i)=L^2(U_1(2r),S)$, 
c.f. \cite[page 257]{Reed&Simon1th-1}, so the above definition shows 
that $\im(\bar D_{2r}+i)=\mathcal H$. Therefore $D_{2r}$, 
as an operator on $H$ is self-adjoint and regular. Consequently we can apply 
the functional calculus of the previous section to define 
the bounded operator $h(D_{2r})$ and $e^{itD_{2r}}$ on $L^2(U_1(2r),W)$.  
The point is that the spectrum of $D_{2r}$ is outside of the 
interval $(-\kappa_0,\kappa_0)$. Since $h$ goes to zero at infinity, 
if $\kappa_0$ is sufficiently 
big then $\|h(D_{2r})\|$ is arbitrary small. 
Let $\sigma$ be a smooth section of $W$ supported in $U_r$. The smooth sections  
$e^{itD}\sigma$ and $e^{itD_{2r}}\sigma$ both satisfy the same wave equation with the same 
initial condition provided $t$ be smaller than $r$. Here we have used the unit 
speed propagation property of the theorem \ref{fpsp}. The uniqueness of the wave 
operator implies their equality for $0\leq t\leq r$. 
Now using the relation \eqref{bradar1} we conclude 
the equality $\phi_1 h(D)\phi_r=\phi_1 h(D_{2r})\phi_r$  
which implies the invertibility of the operator \eqref{jagozari} 
and the vanishing of its index in $K_0(A)$. 
\end{pf}

Now let $M$ be a spin manifold and $N$ be an enlargeable partitioning hypersurface 
of $M$. Moreover assume that there is a smooth map $\phi:M\to N$ such that 
its restriction to $N$ is of non-zero degree.  
Under this condition there is a flat Hilbert $A$-module bundle $V$ on $N$ 
with the following properties (see \cite{HaSc1, HaSc2} and more explicitly 
\cite[theorem 3.1]{Zadeh7}) :
\begin{enumerate}
\item the index of the spin Dirac operator of $N$ twisted by $V$ is a 
non-zero element of $K_0(A)$, 
\item the index of the spin Dirac operator of $N$ twisted by 
$\phi_{|N}^*\,V$ is equal to the index of the spin Dirac operator 
multiplied by $\deg \phi_{|N}$. 
\end{enumerate}
Now by applying the formula \eqref{asli} to the Clifford bundle 
$W:=S\otimes\phi^* V$ we conclude the non vanishing of $\ind(D,N)$. 
This fact along the above theorem show that the scalar curvature of $g$ 
cannot be uniformly positive outside a compact subset of $M$. So we get the 
following theorem 
\begin{thm}
Let $(M,g)$ be a complete riemannian spin manifold and let $N\subset M$ be an 
area-enlargeable partitioning hypersurface. If there is a smooth map 
$\phi\colon M\to N$ such that its restriction to $N$ is of non-zero degree 
then the scalar curvature of $g$ cannot be uniformly positive outside a compact 
subset of $M$.  
\end{thm}
The following corollary is a direct consequence of this theorem
\begin{thm}\label{impthm}
Let $(M,g)$ be a non-compact orientable complete spin $n$-manifold. 
Let $N$ be a $(n-1)$-dimensional sub-manifold of $M$ which is
area-enlargeable. Let $M-N=M_+\sqcup M_-$ and, say, $M_+$ is not compact. 
If there is a map $\phi\colon \overline M_+\to
N$ such that its restriction to $N$ has non-zero degree, then the
scalar curvature of $g$ cannot be uniformly positive 
outside a compact subset of $M_+$.
\end{thm}
\begin{pf}
 By deforming the riemannian metric $g$ in a compact collar neighborhood 
$N\times [0,1)$ in $\overline M_+$ we can and will assume that $g$ takes the 
product form $g_N+dt^2$ where $g_N$ is a riemannian metric on $N$. 
Similarly we can deform $\phi$ in the same collar neighborhood and assume its 
restriction to $N\times[0,1)$ is independent of $t$. 
Let $M_+^-$ denote the (non-complete) 
riemannian manifold $M_+$ with reversed orientation. The riemannian metric $g_{|M_+}$ 
extends naturally (by reflection) to a complete riemannian  metric $g'$ on the manifold 
$M_+\sqcup M_+^-$. The scalar curvature of $g'$ is uniformly positive at infinity 
provided that the scalar curvature of $g$ is uniformly positive at infinity. 
Moreover the map $\phi$ extend (by reflection) to a smooth map $\tilde\phi:=\phi\sqcup\phi$ from 
$M_+\sqcup M_+^-$ onto $N$ with $\deg \tilde\phi_{|N}\neq0$. 
Now we can apply the above theorem to deduce that the scalar curvature of $g$ can not 
be uniformly positive outside a compact. 
\end{pf}

The set $M_+$ satisfying the above properties is called a {\it bad-end} for $M$.  
M. Gromov and B. Lawson proved the 
assertion of the above theorem under the additional condition that 
the Ricci curvature of $g$ be bounded from below on $M_+$, c.f \cite[Theorem 7.46]{GrLa3}. 
Therefore the above theorem is a considerable improvement of their result. 
Moreover, the above theorem improves the following theorem 
of Gromov-Lawson, c.f. \cite[theorem 7.44]{GrLa3} 
(for if $N_0$ is area-enlargeable then $N_0\times S^1$ is area-enlargeable too). 

\emph{Let $(M,g)$ be a connected complete riemannian  manifold which contain a 
compact hypersurface 
$N$ such that: 
$N$ is diffeomorphic to $N_0\times S^1$ where $N_0$ is enlargeable, 
$\pi_1(N_0)\to\pi_1(M)$ is injective and 
there is a non-compact component $M_+$ of $M-N$ and a map $\overline{M}_+\to N$ 
such that its restriction to $N$ has non-zero degree. If one of the following two 
conditions is satisfied then the scalar curvature of $g$ cannot be uniformly 
positive on $M$ 
\begin{enumerate}
\item The map $\overline M_+\to N$ is bounded,
\item $N_0$ has no transversal of finite area.  
\end{enumerate}}
Here a transversal to $N_0$ is a properly embedded 2-dimensional sub manifold of  
$M$ which is transversal to $N_0$ with non-zero intersection number.
By comparing this theorem with the theorem \ref{impthm} it is clear that we have 
relaxed the strong condition on the topology of $N$ (it has not to be of 
the product form $N_0\times S^1$) and we need no condition on the fundamental groups 
in our theorem. Moreover we have relax the boundedness condition of the map 
$\overline M_+\to N$. In addition we have the stronger result that the scalar 
curvature cannot be positive even at infinity. Of course we have paid a price for 
all these: we have assumed $M$ be a spin manifold while in the Gromov-Lawson theorem 
the spin condition is implicit in the enlargeability condition on $N$.

\end{document}